\documentclass[11pt,twoside,reqno]{amsart}

\usepackage{amsthm, amsmath, amscd, amssymb,centernot}
\usepackage[left=2cm,top=2.5cm,bottom=3cm,right=2cm]{geometry}
\usepackage{thmtools}

\usepackage{mathtools}
\DeclarePairedDelimiter{\norm}{\lVert}{\rVert}
\usepackage{tikz-cd}

\usepackage{xcolor} 
\colorlet{mdtRed}{red!50!black}
\definecolor{dblue}{rgb}{0,0,.6}
\usepackage[colorlinks]{hyperref}
\hypersetup{linkcolor=blue,citecolor=dblue,filecolor=dullmagenta,urlcolor=mdtRed}
\usepackage[capitalise]{cleveref}
\newcommand{\mb}{\mathbb}

\newcommand{\Xr}{\widetilde{X}_r(x)}

\newcommand{\ex}{\varepsilon_x}
\newcommand{\mx}{{\rm mult}_x}
\theoremstyle{plain}
\numberwithin{equation}{section}
\newtheorem{theorem}{Theorem}[section]
\newtheorem*{theorem*}{Theorem}
\newtheorem{proposition}[theorem]{Proposition}
\newtheorem{lemma}[theorem]{Lemma}
\newtheorem{corollary}[theorem]{Corollary}
\newtheorem{conjecture}[theorem]{Conjecture}
\newtheorem{definition}[theorem]{Definition}
\theoremstyle{definition}

\title{Seshadri Functions on Surfaces}
\author{Cyril J. Jacob} 
	\address{Department of Mathematics, Indian Institute of Technology Bombay, Powai, Mumbai 
		400076, Maharashtra, India}
	\email{cyriljjacob@iitb.ac.in}

\begin{document}
\date{\today}

	\begin{abstract}
	    We study the Seshadri function on the ample cone of a smooth complex projective surface at a fixed point. At a submaximal class, we show that the function is differentiable if and only if there is a unique Seshadri curve. At a maximal class $L$, differentiability is equivalent to every curve $C$ computing the Seshadri constant satisfying $\frac{[C]}{{\rm mult}_x C}=\frac{[L]}{\sqrt{L^{2}}}$. We also prove that the non-differentiability locus is nowhere dense.
        
        For the blow-up of $\mathbb P^{2}$ at $r\ge 9$ very general points, assuming the $(-1)$-curves conjecture for the blow-up at $r+1$ very general points, we construct a non-empty open subset of the ample cone on which the Seshadri function at a very general point equals $\sqrt{L^{2}}$. Consequently, this function is not piecewise linear. Under the assumption of the Strong SHGH conjecture, we obtain a common non-empty open subset on which maximality holds at every point when $r\ge36$.
	\end{abstract}

	\maketitle
    \section{Introduction}
    Throughout the article, we work over the field of complex numbers. By a \textit{surface}, we mean a smooth integral projective variety of dimension two over $\mathbb C$.

    Seshadri constants are numerical invariants associated to a nef line bundle on a projective surface. We start by recalling the definition of the Seshadri constant of a nef line bundle at a point.
    \begin{definition}
		Consider a smooth complex projective surface $X$. The \textit{Seshadri constant} of a nef line bundle $L$ on $X$ at a point $x \in X$ is defined as follows:
		$$\varepsilon(X,L,x)=\inf\bigg \{ \frac{L\cdot  C}{\text{mult}_x C}~~\bigg  |  ~~C \text{ is a reduced irreducible curve in $X$ passing through $x$} \bigg  \}. $$
	\end{definition}
    This definition extends to real nef
numerical classes. It depends only on the numerical class of $L$ and satisfies
$0\leq \varepsilon(X,L,x)\leq \sqrt{L^2}.$
An ample class $L$ is called \emph{submaximal at $x$} if the upper
bound is strict, and \emph{maximal at $x$} if equality holds.
A curve attaining the infimum is called a \emph{Seshadri curve}
for $L$ at $x$. More details about Seshadri constants and local positivity can be found in \cite{Prim09}.

    	The study of Seshadri constants has led to the introduction of an associated function.
	Fixing a point $x\in X$ and varying $L$
gives the \emph{Seshadri function}
$$\varepsilon_x:{\rm Nef}_{\mathbb R}(X)\longrightarrow
\mathbb R,\qquad L\longmapsto\varepsilon(X,L,x),$$
	where $\text{\rm Nef}_{\mathbb{R}}(X)$ denotes the nef cone of $X$ over $\mathbb{R}$.
    
    Each curve through $x$ contributes a linear function
$L\mapsto L\cdot C/{\rm mult}_x C$, and $\varepsilon_x$
is their infimum. This description suggests that changes in the
curves computing the Seshadri constant should be reflected in the
regularity of $\varepsilon_x$. At maximal classes, however, the intersection form enters through the upper bound
$\sqrt{L^2}$. The Seshadri function therefore brings together
linear contributions from curves and the geometry of this bound.

Seshadri functions on abelian surfaces have been studied in
\cite{BS08,BSc22,Sch}. On any smooth projective surface,
$\varepsilon_x$ is continuous on the nef cone and is locally
piecewise linear near every submaximal ample class; see
\cite[Section~1]{Sch}. For an abelian surface of Picard rank
at least two, M. Schmidt also proved that the maximality locus
$$W_x(X)=\left\{L\in{\rm Nef}_{\mathbb R}(X)\;\middle|\;\varepsilon(X,L,x)=\sqrt{L^2}\right\}$$
is nowhere dense \cite[Lemma~4.1]{Sch}.

These results raise two related questions. Can the differentiability
of the Seshadri function be characterized by the curves computing
its value? And can the maximality locus contain an open set, so
that the function agrees with $\sqrt{L^2}$ throughout a region
of the ample cone? We address the first question for arbitrary
smooth projective surfaces and the second for blow-ups of the
plane under suitable SHGH assumptions.

    We begin by studying differentiability on the ample cone. In \cref{sec:differential}, we first prove that $\varepsilon_x$ is concave and
locally Lipschitz continuous on ${\rm Amp}_{\mathbb R}(X)$.
Rademacher's theorem then implies that it is differentiable almost
everywhere. We give precise criteria for differentiability at both
submaximal and maximal ample classes. At a submaximal class $L$,
the function $\varepsilon_x$ is differentiable if and only if
$L$ has a unique Seshadri curve at $x$. At a maximal class $L$,
it is differentiable if and only if every curve $C$ computing
$\varepsilon_x(L)$ satisfies
$$
\frac{[C]}{{\rm mult}_x C}
=
\frac{[L]}{\sqrt{L^2}}.
$$
In particular, $\varepsilon_x$ is differentiable at a maximal
class whose Seshadri constant is not computed by any curve.

    Then we study the points $L\in{\rm Amp}_{\mb R}(X)$ where $\ex$ is not differentiable. We define
    $$ N_x:=\left\{ L \in {\rm Amp}_{\mb R}(X) ~\middle|~ \ex \text{ is not differentiable at } L \right\}.$$

    Although Rademacher’s theorem shows that this set has measure zero, that conclusion alone does not exclude density. We also establish the following topological property.

    \begin{theorem} (\cref{thm: N_x nowhere dense})
    For a smooth surface $X$, 
        the set $N_x$ is nowhere dense in ${\rm Amp}_{\mb R}(X)$ for all $x\in X$.
    \end{theorem}

    In \cref{sec:Submaximality}, we study maximality on rational
surfaces. Let $X_r$ denote the blow-up of $\mathbb P^2$ at $r$
very general points. In contrast with the abelian-surface
result above, we show that the maximality locus on $X_r$ can
have nonempty interior under the Weak SHGH conjecture.
	\begin{theorem}(\cref{open})
		For $r\geq9$, assume that Conjecture~2.9 holds for $X_{r+1}$.
Then, for a very general point $x\in X_r$, there exists a
nonempty open subset
$U_x\subset{\rm Amp}_{\mathbb R}(X_r)$ such that
$$\varepsilon(X_r,L,x)=\sqrt{L^2}
\qquad\text{for every }L\in U_x.$$
In particular, the Seshadri function $\varepsilon_x$ is not
piecewise linear. Moreover, $U_x$ contains the numerical class
of an integral ample line bundle whose Seshadri constant at
$x$ is irrational.
	\end{theorem}
	
	Our construction also gives infinitely many integral ample
line bundles with pairwise non-proportional numerical classes
and irrational Seshadri constants at a very general point.
The existence of infinitely many ample line bundles with irrational
Seshadri constants under the Weak SHGH conjecture was already
observed in \cite[Remark~2.5]{Han1}; here we give an explicit construction
with pairwise non-proportional numerical classes.

Finally, under the Strong SHGH conjecture for $X_r$, we prove
that for $r\geq36$ there exists a common nonempty open subset
$U\subset{\rm Amp}_{\mathbb R}(X_r)$ such that
$$\varepsilon(X_r,L,x)=\sqrt{L^2}
\qquad
\text{for every }L\in U\text{ and every }x\in X_r.$$
Thus the same open region consists of classes that are maximal
at every point of the surface.

The paper is organized as follows. \cref{sec:prelim} recalls the
preliminary results and conjectures used in the proofs.
\cref{sec:differential} establishes the differentiability criteria and
the nowhere density of the non-differentiability locus.
\cref{sec:Submaximality} constructs open submaximality and maximality regions
on blow-ups of $\mb P^2$ and proves the results on irrational
Seshadri constants.
    \subsection*{Acknowledgments}
	   The author is grateful to Prof. Krishna Hanumanthu for suggesting the problem and for many helpful discussions, suggestions, and comments. 

       \subsection*{ AI Disclosure} The author acknowledges the use of ChatGPT for assistance with literature discovery, the exploration of ideas, language editing, and final proofreading of the manuscript prior to posting. The author assumes full responsibility for the mathematical content, arguments, and claims presented in this work.

\section{Preliminaries}\label{sec:prelim}

In this section, we recall the blow-up characterization of Seshadri constants, the Weyl group action on plane blow-ups, and the SHGH-type conjectures used in our constructions. The first part relates submaximal Seshadri constants to negative curves. The remaining parts introduce the tools used to control these curves on blow-ups of $\mb P^2$.

\subsection{Seshadri constants and negative curves}

We begin with the following characterization of Seshadri constants.
\begin{lemma}\label{alt}
	Let $L$ be a nef line bundle on a surface $X$. Consider the blow-up $\pi_x:\widetilde{X}(x)\rightarrow X$ of $X$ at a point $x\in X$, and let $E_x=\pi_x^{-1}(\{x\})$ denote the exceptional divisor. Then
	$$\varepsilon(X,L,x)=\sup\left\{t\in\mathbb{R}_{\ge 0}~~\Big|~~\pi_x^{\ast}L-tE_x \text{ is nef on } \widetilde{X}(x)\right\}.$$
\end{lemma}

Since a nef class on a surface has non-negative self-intersection, Lemma \ref{alt} gives
$$\varepsilon(X,L,x)\leq\sqrt{L^2}.$$
We recall the terminology associated with this upper bound.
\begin{definition}
	Let $L$ be an ample line bundle on a smooth projective surface $X$ and let $x\in X$. We say that \textit{$L$ is submaximal at $x$} if
	$$\varepsilon(X,L,x)<\sqrt{L^2},$$
	and \textit{maximal at $x$} if
	$$\varepsilon(X,L,x)=\sqrt{L^2}.$$
\end{definition}

If $L$ is submaximal at $x$, then there exists an integral curve $C$ passing through $x$ such that
$\varepsilon(X,L,x)=\frac{L\cdot C}{\text{mult}_xC}$
(see \cite[Proposition 5.1]{TB99}, \cite[Lemma 1.9]{TS01}, \cite[Corollary 2]{KO02}). Such a curve is called a \textit{Seshadri curve} of $L$ at $x$. In the submaximal case, T. Szemberg proved that the number of such curves is at most the Picard rank $\rho(X)$ \cite[Proposition 1.7]{TS01}.

Suppose that $L$ is submaximal at $x$ and that $C$ is a Seshadri curve. If $\widetilde{C}$ is its strict transform on $\widetilde{X}(x)$, then the Hodge index theorem gives $\widetilde{C}^2<0$. 

For a positive integer $m$, an integral curve $C$ on $X$ is called a \textit{$(-m)$-curve} if
$$C^2=-m\qquad\text{and}\qquad K_X\cdot C=m-2,$$
where $K_X$ is the canonical line bundle of $X$.

\subsection{Plane blow-ups and the Weyl group}

Let $\pi:X_r\rightarrow\mb P^2$ be the blow-up of $\mb P^2$ at $r$ very general points $p_1,p_2,\dots,p_r$. Let $H$ be the pullback of a line avoiding these points, and let $E_i=\pi^{-1}(\{p_i\})$ denote the exceptional divisor corresponding to $p_i$, for $1\leq i\leq r$. The linear equivalence classes of $H,E_1,E_2,\dots,E_r$ form a basis for $\text{Pic}(X_r)$, with intersection numbers
$$H^2=1,\qquad E_i^2=-1,\qquad H\cdot E_i=0,\qquad E_i\cdot E_j=0\text{ for }i\neq j.$$
In this basis, the canonical class is
$$K_{X_r}=-3H+E_1+E_2+\dots+E_r.$$

We write
$$\text{Pic}_{\mb R}(X_r)=\text{Pic}(X_r)\otimes_{\mb Z}\mb R,\qquad
\text{Pic}_{\mb Q}(X_r)=\text{Pic}(X_r)\otimes_{\mb Z}\mb Q.$$
On $X_r$, these spaces are naturally identified with the corresponding spaces of numerical divisor classes.

We will also consider the blow-up $\pi_x:\Xr\to X_r$ at a point $x\in X_r$. We retain the notation $H,E_1,\dots,E_r$ for their pullbacks to $\Xr$ and write $E_x:=\pi_x^{-1}(\{x\})$ for the new exceptional divisor.

Assume $r\geq3$. The Weyl group action allows us to arrange divisor classes into standard forms while preserving their intersection numbers. Consider the simple roots
$$\Gamma_0=H-E_1-E_2-E_3,\qquad
\Gamma_i=E_i-E_{i+1}\text{ for }1\leq i\leq r-1.$$
The associated Weyl group $W_r$ is the subgroup of ${\rm Aut}({\rm Pic}(X_r))$ generated by the reflections
$$\sigma_i:\text{Pic}(X_r)\to\text{Pic}(X_r),\qquad
\sigma_i(D)=D+(D\cdot\Gamma_i)\Gamma_i,$$
for $0\leq i\leq r-1$.

Every element of $W_r$ corresponds to a Cremona transformation on $\mb P^2$ (see \cite[Theorem 1]{Dol}). We will use the following properties of this action.
\begin{lemma}\cite[Lemma 1.5]{Nag}, \cite[Theorem 1]{Dol}\label{mainlemma}
	For all $\sigma\in W_r$ and $C,D\in\text{Pic}(X_r)$, we have:
	\begin{enumerate}
		\item $C\cdot D=\sigma(C)\cdot\sigma(D)$;
		\item $\sigma(-K_{X_r})=-K_{X_r}$;
		\item $h^0(D)=h^0(\sigma(D))$;
		\item if $D$ is represented by an irreducible curve, then so is $\sigma(D)$.
	\end{enumerate}
\end{lemma}

The Weyl group action extends linearly to $\text{Pic}_{\mb R}(X_r)$. The following standard forms will be used to control intersections with $(-1)$-curves in our constructions of maximality regions.
\begin{definition}
	Assume $r\geq3$. A line bundle $L=dH-m_1E_1-m_2E_2-\dots-m_rE_r$ on $X_r$ is called \textit{pseudo-standard} if
	$$d\geq0,\qquad m_1\geq m_2\geq\dots\geq m_r,\qquad d\geq m_1+m_2+m_3.$$
	It is called \textit{standard} if it is pseudo-standard and $m_r\geq0$.
\end{definition}

The following two propositions explain how effective and nef divisor classes can be expressed in these forms. We include proofs of both propositions for completeness, even though the statements are well known.
\begin{proposition}\label{stdcomes}
	Assume $r\geq3$. Let $D=dH-m_1E_1-\dots-m_rE_r$ be a divisor on $X_r$ such that $m_1\geq m_2\geq\dots\geq m_r$. Then there exists $\sigma\in W_r$ such that either $\sigma(D)$ is pseudo-standard or $\sigma(D)\cdot H<0$. In particular, if $D$ is effective, there exists $\sigma\in W_r$ such that $\sigma(D)$ is in pseudo-standard form.
\end{proposition}

\begin{proof}
	If $d<0$ or $D$ is already in pseudo-standard form, there is nothing to prove. Otherwise, $d\geq0$ and $d<m_1+m_2+m_3$. Set
	$$d'=2d-m_1-m_2-m_3,$$
	$$m_i'=d-m_1-m_2-m_3+m_i\text{ for }i=1,2,3,\qquad
	m_i'=m_i\text{ for }i>3.$$
	Then
	$$\sigma_0(D)=d'H-\sum_{i=1}^{r}m_i'E_i.$$
	The reflections $\sigma_i$ for $1\leq i\leq r-1$ interchange adjacent multiplicity coefficients, so we can reorder the new coefficients in decreasing order. Moreover,
	$$d'=d+(d-m_1-m_2-m_3)<d.$$
	Thus, whenever the resulting class has non-negative degree and is not pseudo-standard, we can repeat the operation and strictly decrease its degree. Since the degree is an integer, this process terminates after finitely many steps, yielding either a pseudo-standard divisor or a divisor of negative degree.

	If $D$ is effective, Lemma \ref{mainlemma} implies that every class in its Weyl group orbit is effective. Such a class has non-negative intersection with $H$, so the negative-degree alternative cannot occur.
\end{proof}

\begin{proposition}\label{nefstd}
	Assume $r\geq3$. Let $D=dH-m_1E_1-\dots-m_rE_r$ be a nef divisor on $X_r$. Then there exists $\sigma\in W_r$ such that $\sigma(D)$ is standard.
\end{proposition}

\begin{proof}
	For every $\sigma\in W_r$, we have
	$\sigma(D)\cdot H=D\cdot\sigma^{-1}(H)\geq0,$
	since $D$ is nef and $\sigma^{-1}(H)$ is effective. After reordering the coefficients of $D$ if necessary, Proposition \ref{stdcomes} therefore gives an element $\sigma\in W_r$ such that $\sigma(D)$ is pseudo-standard.

	Moreover,
	$\sigma(D)\cdot E_r=D\cdot\sigma^{-1}(E_r)\geq0,$
	because $\sigma^{-1}(E_r)$ is effective. Thus the smallest multiplicity coefficient of $\sigma(D)$ is non-negative, and $\sigma(D)$ is standard.
\end{proof}

\subsection{SHGH-type conjectures}

To establish maximality of Seshadri constants, we will use conjectural restrictions on negative curves on plane blow-ups and on curves with a point of prescribed multiplicity. These restrictions arise from the interpolation problem for plane curves: determining the dimension of a linear system with assigned multiplicities at very general points.

Consider the linear system
$\mathcal{L}=|dH-m_1E_1-m_2E_2-\dots-m_rE_r|$
on $X_r$, where $d\geq1$ and $m_1,m_2,\dots,m_r\geq0$ are integers. Counting the conditions imposed by the prescribed multiplicities gives the \textit{expected dimension}
$$\text{edim }\mathcal{L}
=\max\Bigg\{{d+2\choose2}
-\displaystyle\sum_{i=1}^{r}{m_i+1\choose2}-1,-1\Bigg\}.$$
We say that $\mathcal{L}$ is \textit{non-special} if $\dim\mathcal{L}=\text{edim }\mathcal{L}$ and \textit{special} otherwise.

The SHGH conjecture describes special linear systems on blow-ups of the plane at very general points. It was independently proposed by Segre \cite{Ser}, Harbourne \cite{Har1}, Gimigliano \cite{Gim} and Hirschowitz \cite{Hir}. We use the following formulation due to Segre.
\begin{conjecture}[SHGH Conjecture]\label{SHGH}
	If $\mathcal{L}$ is special, then every divisor in $\mathcal{L}$ is non-reduced.
\end{conjecture}

The SHGH conjecture is known for $r\leq9$. To study Seshadri constants at arbitrary points, we also use the following stronger conjecture, which controls curves having a point of prescribed multiplicity. It is called the \textit{Strong SHGH conjecture} or the \textit{Generalized SHGH conjecture} \cite[Conjecture 3.6]{Far}.
\begin{conjecture}[Generalized SHGH or Strong SHGH]\label{SSHGH}
	Let $X_r$ be the blow-up of $\mb P^2$ at $r$ very general points, and let $d\geq1$, $t\geq1$ and $m_1,m_2,\dots,m_r\geq0$ be integers such that
	$${d+2\choose2}
	-\displaystyle\sum_{i=1}^{r}{m_i+1\choose2}
	\leq\max\Bigg\{{t+1\choose2}-2,0\Bigg\}.$$
	Then any curve $C\in|dH-m_1E_1-\dots-m_rE_r|$ having a point of multiplicity $t$ is non-reduced.
\end{conjecture}

The case $t=1$ is equivalent to the SHGH conjecture. We will also use the following consequence of SHGH, known as the \textit{Weak SHGH conjecture} or the \textit{$(-1)$-curves conjecture}.
\begin{conjecture}[Weak SHGH]\label{WSHGH}
	Let $X_r$ be the blow-up of $\mb P^2$ at $r$ very general points. If $C$ is an irreducible and reduced curve on $X_r$ with $C^2<0$, then $C$ is a $(-1)$-curve.
\end{conjecture}

Conjectures \ref{SHGH} and \ref{WSHGH} remain open for $r\geq 10$,
while Conjecture \ref{SSHGH} remains open even for $r=1$.
These conjectures are related to Nagata's conjecture by the implications
\[
\cref{SSHGH}\implies \cref{SHGH}\implies \cref{WSHGH}
\implies \text{Nagata's conjecture}.
\]
Nagata's conjecture states that the class
$\sqrt{r}H-\sum_{i=1}^r E_i$ is nef on the blow-up of $\mathbb{P}^2$
at $r\geq 10$ very general points. A proof has recently been
announced by OpenAI; see also Tevelev's exposition \cite{T26}.
Our proofs of maximality also require control of negative curves
on a further blow-up, so we retain the Weak or Strong SHGH
hypotheses in the corresponding statements.


    \section{General properties of $\ex$}\label{sec:differential}
    Throughout this section, we work with ample real numerical classes unless integrality is explicitly stated. We first express $\varepsilon_x$ as the infimum of linear functionals on the ample cone.

    Observe that
    $$\varepsilon_x(L)=\inf_C \frac{L\cdot C}{{\rm mult}_x C},$$
    where the infimum varies over integral curves $C$ passing through $x$. Note that $\frac{\_\cdot C}{\mx C}$ is a linear functional on $N_1(X)_{\mb R}$.

    \begin{proposition}
        $\varepsilon_x$ is concave.
    \end{proposition}
    \begin{proof}
        Let $L$ and $M$ be ample line bundles on $X$. For $t\in [0,1]$, we have
        \begin{equation*}
            \begin{split}
                \varepsilon_x(tL+(1-t)M)&=\inf_C\frac{(tL+(1-t)M)\cdot C}{{\rm mult}_x C}\\
                &\ge t\inf_C\frac{L\cdot C}{{\rm mult}_x C}+(1-t)\inf_C\frac{M\cdot C}{{\rm mult}_x C}\\
                &=t\varepsilon_x (L)+(1-t)\varepsilon_x(M).
            \end{split}
        \end{equation*}
    \end{proof}
    \begin{definition}
        A function $f: U \subseteq \mb R^n \to \mb R$ is called \textit{locally Lipschitz} if for each $p\in U$, there exists an open subset $V$ of $U$ and a constant $K$ such that $p\in V$ and 
        $$|f(u)-f(v)|\le K\|u-v\| \hspace{.5 cm} \text{for all }u,v\in V.$$
    \end{definition}

    \begin{proposition}\label{prop:lipschitz}
    Let $X$ be a smooth surface and $x\in X$ be a point. Then the Seshadri function $\ex$ at $x$ is locally Lipschitz continuous on ${\rm Amp}_{\mb R}(X)$.
    \end{proposition}
    \begin{proof}
    Consider an ample line bundle $L$ on $X$.
    Choose $\gamma>0$ such that $\overline{B_\gamma(L)}\subset {\rm Amp}_{\mb R}(X)$. Let $K$ denote the maximum value attained by the continuous function $\ex$ on $\overline{B_\gamma(L)}$. Choose $0<\delta<\gamma/2$. Consider distinct $M_1,M_2\in B_\delta(L).$ 
    Define $$M=M_2+\frac{\delta}{\lVert M_2-M_1\rVert}(M_2-M_1).$$
    Then we have
    $$\norm{M-L}\le \norm{M_2-L}+ \frac{\delta}{\lVert M_2-M_1\rVert} \norm{M_2-M_1}<2\delta<\gamma.$$
    As $\overline{B_\gamma(L)}\subset {\rm Amp}_{\mb R}(X)$, we conclude that $M$ is ample. 
    
    Further,
     $${\lVert M_2-M_1\rVert}M={\lVert M_2-M_1\rVert}M_2+{\delta}(M_2-M_1)=(\lVert M_2-M_1\rVert+\delta)M_2-{\delta}M_1.$$
     Hence,
    $$M_2=\frac{\delta}{\delta+\norm{M_2-M_1}}M_1+\frac{\norm{M_2-M_1}}{\delta+\norm{M_2-M_1}}M.$$

    Therefore from the concavity of the function $\ex$, we have:
        $$\ex(M_2)\ge \frac{\delta}{\delta+\norm{M_2-M_1}}\ex(M_1)+
        \frac{\norm{M_2-M_1}}{\delta+\norm{M_2-M_1}}\ex(M)\ge \frac{\delta}{\delta+\norm{M_2-M_1}}\ex(M_1).$$
        The last inequality follows from $\ex(M)> 0$.
        Rearranging gives us:
        $$\ex(M_1)\le \frac{\delta+\norm{M_2-M_1}}{\delta}\ex(M_2) = \ex(M_2)+ \frac{\norm{M_2-M_1}}{\delta}\ex(M_2).$$
        Thus,
        $$\ex(M_1)- \ex(M_2)\le \frac{\ex(M_2)}{\delta}\norm{M_2-M_1}\le \frac{K}{\delta}\norm{M_2-M_1}. $$

        Now, if we interchange the roles of $M_1$ and $M_2$, we get
         $$\ex(M_2)- \ex(M_1)\le \frac{\ex(M_1)}{\delta}\norm{M_1-M_2}\le \frac{K}{\delta}\norm{M_2-M_1} .$$
         Hence,
         $$|\ex(M_1)- \ex(M_2)|\le \frac{K}{\delta}\norm{M_2-M_1} .$$

    \end{proof}

    \begin{theorem}(Rademacher's Theorem)~~\cite[Theorem 3.3]{EG15}\label{thm:Rademacher}
        Let $f:U \subset \mb R^n \to \mb R$ be a locally Lipschitz continuous function where $U$ is an open subset of $\mb R^n$. Then $f$ is differentiable outside a set of zero Lebesgue measure.
    \end{theorem}

    \begin{proposition}\label{prop:lebesgue measure}
        The function $\ex$ is differentiable on ${\rm Amp}_{\mb R}(X)$ outside a set of zero Lebesgue measure.
    \end{proposition}
    \begin{proof}
        This follows from Theorem \ref{thm:Rademacher} and  \cref{prop:lipschitz}.
    \end{proof}

    For further study, let 
    $$ N_x:=\left\{ L \in {\rm Amp}_{\mb R}(X) \middle| \ex \text{ is not differentiable at } L \right\}.$$
    From Proposition \ref{prop:lebesgue measure}, we know that $N_x$ has Lebesgue measure zero. Our goal is to understand the set $N_x$ better.

    \begin{proposition}
    \label{prop:diffrentiablity in  submaximal case}
        Let $X$ be a smooth surface and $L$ be an ample line bundle on $X$ such that $\ex(L)<\sqrt{L^2}$ at a point $x\in X$. Let $C_1,\dots, C_k$ be all Seshadri curves of $L$ at $x$. Then $\ex$ is differentiable at $L$ if and only if 
        $$\frac{[C_1]}{{\rm mult}_x C_1}=\frac{[C_2]}{{\rm mult}_x C_2}=\dots = \frac{[C_k]}{{\rm mult}_x C_k} \in N_1(X)_{\mb R}. $$
        Further,
 $$d\varepsilon_x|_L(M)=\frac{M\cdot C_i}{{\rm mult}_x C_i} \hspace{.5 cm}
\text{for every }M\in N_1(X)_{\mathbb R}
\text{ and every }i.$$
        
    \end{proposition}
    \begin{proof}
        As $\ex(L)<\sqrt{L^2}$, by \cite[Proposition 1.3]{Sch} we can see that there exist integral curves $C_1,C_2,\dots,C_k$ that are Seshadri curves of $L$ at $x$ and a neighborhood $U\subset {\rm Nef}_{\mb R}(X)$ of $L$ such that $$\ex(L')=\min_{1\le i\le k} \frac{L'\cdot C_i}{{\rm mult}_x C_i}$$ for all $L'\in U$.

        As $U$ is an open set, for each $M\in N_1(X)_{\mb R}$, for sufficiently small $t>0$ we have $L+tM\in U$ .
        Since $\ex(L)=\frac{L\cdot C_i}{{\rm mult}_x C_i}$ for all $1\le i\le k$, we have:
        $$\ex(L+tM)=\min_{1\le i\le k} \frac{(L+tM)\cdot C_i}{{\rm mult}_x C_i}=\ex(L)+\min_{1\le i\le k} t\frac{M\cdot C_i}{{\rm mult}_x C_i}=\ex(L)+t\min_{1\le i\le k} \frac{M\cdot C_i}{{\rm mult}_x C_i}.$$
        The last equality is because $t>0$.
        
        Therefore,
               $$\lim_{t \to 0^+}\frac{\ex(L+tM)-\ex(L)}{t} =\min_{1\le i\le k} \frac{M\cdot C_i}{{\rm mult}_x C_i}.$$

        Similarly for $t<0$ with small $|t|$, we get:
        $$\ex(L+tM)=\min_{1\le i\le k} \frac{(L+tM)\cdot C_i}{{\rm mult}_x C_i}=\ex(L)+\min_{1\le i\le k} t\frac{M\cdot C_i}{{\rm mult}_x C_i}=\ex(L)+t\max_{1\le i\le k} \frac{M\cdot C_i}{{\rm mult}_x C_i}.$$
        The last equality follows because $t<0$.

        Therefore,
               $$\lim_{t \to 0^-}\frac{\ex(L+tM)-\ex(L)}{t} =\max_{1\le i\le k} \frac{M\cdot C_i}{{\rm mult}_x C_i}.$$

Suppose first that $\varepsilon_x$ is differentiable at $L$.
Then the two one-sided limits above agree for every
$M \in N^1(X)_{\mathbb R}$. Hence
$$
\min_{1\leq i\leq k}
\frac{M\cdot C_i}{{\rm mult}_x C_i}
=
\max_{1\leq i\leq k}
\frac{M\cdot C_i}{{\rm mult}_x C_i}.
$$
Therefore, for every $i,j$ and every
$M \in N^1(X)_{\mathbb R}$,
$$
M\cdot\left(
\frac{[C_i]}{{\rm mult}_x C_i}
-
\frac{[C_j]}{{\rm mult}_x C_j}
\right)=0.
$$
Therefore we have,
$$
\frac{[C_1]}{{\rm mult}_x C_1}
=\cdots=
\frac{[C_k]}{{\rm mult}_x C_k}
\quad\text{in }N^1(X)_{\mathbb R}.
$$

Conversely, suppose the  numerical classes $\frac{[C_i]}{\mx C_i}$
are the same for all $1\le i \le k$, and denote their common value by $\frac{[C]}{\mx C}$.
The local expression for $\varepsilon_x$ gives
$$
\varepsilon_x(L')
=
\min_{1\leq i\leq k}
\frac{L'\cdot C_i}{{\rm mult}_x C_i}
=
\frac{L'\cdot C}{\mx C}
\qquad\text{for every }L'\in U.
$$
Thus $\varepsilon_x$ agrees with a linear functional
on a neighborhood of $L$, so it is differentiable at $L$.
Moreover,
$$
d\varepsilon_x|_L(M)
=
M\cdot\frac{[C]}{\mx C}
=
\frac{M\cdot C_i}{{\rm mult}_x C_i}
$$
for every $M\in N^1(X)_{\mathbb R}$ and every
$1\leq i\leq k$.
This completes the proof.
    \end{proof}
    \begin{lemma}\label{lem:seshadri curve class different}
    Let $L$ be an ample line bundle such that $\ex(L)<\sqrt{L^2}$.
        Let $C$ and $D$ be two distinct Seshadri curves of $L$ at $x$. Then
        $$\frac{[C]}{{\rm mult}_x C}\ne \frac{[D]}{{\rm mult}_x D} \text{ in } N_1(X)_{\mb R} .$$
    \end{lemma}
    \begin{proof}
        Suppose, to the contrary, that $$\gamma:=\frac{[C]}{{\rm mult}_x C}= \frac{[D]}{{\rm mult}_x D} \text{ in } N_1(X)_{\mb R}.$$
        Since $C$ and $D$ are Seshadri curves, we have $$C^2-({\rm mult}_x C)^2<0.$$
        Thus, $$\gamma^2=\frac{C^2}{({\rm mult}_x C)^2}<1.$$
        But since $C$ and $D$  are distinct integral curves, we have $C\cdot D \ge ({\rm mult}_x C) ({\rm mult}_x D)$. So we conclude
        $$\gamma^2=\frac{C\cdot D}{({\rm mult}_x C) ({\rm mult}_x D)}\ge 1.$$
        This is a contradiction.
    \end{proof}
    \begin{theorem}
        Let $X$ be a smooth surface and $L$ be an ample line bundle on $X$ such that $\ex(L)<\sqrt{L^2}$.  $\ex$ is differentiable at $L$ if and only if $L$ has a unique Seshadri curve $C$.
        Further,
 $$d\varepsilon_x|_L(M)=\frac{M\cdot C}{{\rm mult}_x C}.$$
    \end{theorem}
    \begin{proof}
        It follows from \cref{prop:diffrentiablity in  submaximal case} and \cref{lem:seshadri curve class different}.
    \end{proof}

    \begin{lemma}\label{lem:bdd numerical class}
        Let $L$ be an ample line bundle on a surface $X$. If  $\gamma_n$ is a sequence of effective numerical classes such that $L\cdot \gamma_n$ is bounded, then $\{\gamma_n\}$ is bounded.
    \end{lemma}
    \begin{proof}
        Suppose, to the contrary, that $\gamma_n$ is unbounded. By passing to a subsequence, without loss of generality, we can assume that $\norm{\gamma_n}\to \infty$. Since $\frac{\gamma_n}{\norm{\gamma_n}}$ belongs to the unit sphere in a finite-dimensional space, by again passing to a subsequence we can assume $\frac{\gamma_n}{\norm{\gamma_n}}$ converges to $\gamma$ with $\norm \gamma=1$. 

        By assumption, the sequence $L\cdot\gamma_n$ is bounded. This enables us to conclude
        $$L\cdot \gamma=\lim_{n \to \infty} \frac{L\cdot \gamma_n}{\norm{\gamma_n}}= 0.$$

        We will show that $\gamma=0$ in $N_1(X)_{\mb R}$, which gives a contradiction as $\norm \gamma =1$. 
        Let $D\in N_1(X)_{\mb R}$. Since $L$ is ample, both $L+tD$ and $L-tD$ are ample for sufficiently small $t>0$. Since $\gamma_n$ are effective, 
        we have $(L\pm tD)\cdot \gamma_n\ge0$. 
 Dividing by $\norm{\gamma_n}$ and taking limits, we have $\pm tD\cdot \gamma\ge 0$. 
Thus we can conclude $D\cdot \gamma =0$ for all $D\in N_1(X)_{\mb R}$. Thus $\gamma=0$ in $N_1(X)_{\mb R}$.
    \end{proof}

    \begin{lemma}\label{lem:seq argument}
        Let $L$ be an ample line bundle on a surface $X$ with $\ex(L)=\sqrt{L^2}$. Suppose every curve $C$ computing $\ex(L)$ satisfies $$\frac{[C]}{{\rm mult}_x C}=\frac{[L]}{\sqrt{L^2}}.$$
        
        Let $C_n$ be a sequence of integral curves passing through $x$ such that $$\frac{L\cdot C_n}{{\rm mult}_x C_n} \to \sqrt{L^2}.$$
        Then $\frac{[C_n]}{\mx C_n}$ converges to $\frac{[L]}{\sqrt{L^2}}.$
    \end{lemma}
    \begin{proof}
        Let $\gamma_n$ denote the numerical class $\frac{[C_n]}{{\rm mult}_x C_n}$ in $N_1(X)_{\mb R}$. Since $L\cdot \gamma_n$ converges, the sequence $L\cdot \gamma_n$ is bounded, and we know that $\gamma_n$ are effective numerical classes. Hence, by \cref{lem:bdd numerical class}, the sequence $(\gamma_n)$ is bounded in $N_1(X)_{\mb R}$. Hence $(\gamma_n)$ has a convergent subsequence. So, without loss of generality, assume that $\gamma_n \to \gamma \in N_1(X)_{\mb R}$. 
        

        Suppose that some class occurs infinitely many times in the sequence
$(\gamma_n)$. Since $\gamma_n\to\gamma$, this class must be $\gamma$.
Choose an index $n$ for which $\gamma_n=\gamma$, and set $C=C_n$.
Then $\gamma=\frac{[C]}{\mx C}$. Since
$L\cdot\gamma=\sqrt{L^2}$, the curve $C$ computes $\varepsilon_x(L)$.
Then by hypothesis, $\gamma=\frac{[L]}{\sqrt{L^2}}$.
    
        Otherwise, without loss of generality, we can suppose  $\gamma_n$ are mutually distinct. In this case, we also show that $\gamma=\frac{[L]}{\sqrt{L^2}}$.
        
        Since $C_n$ are integral and distinct,   for $i\neq j$ we have $$C_i\cdot C_j\ge ({\rm mult}_x C_i){(\rm mult}_x C_j).$$ Thus $\gamma_i\cdot \gamma_j\ge 1$. Hence by taking limits we can see $\gamma^2\ge 1$. Now we can use the Hodge index theorem to see that 
        $$1\le \gamma^2= \frac{L^2\gamma^2}{L^2}\le \frac{(L\cdot \gamma)^2 }{L^2}=1.$$
        The last equality follows from $L\cdot \gamma_n \to L\cdot \gamma= \ex(L)=\sqrt{L^2}$.

        Assume $\gamma\ne\frac{[L]}{\sqrt{L^2}}$. 
        $$L\cdot \left(\gamma-\frac{L}{\sqrt{L^2}}\right)=L\cdot \gamma-\frac{L^2}{\sqrt{L^2}}=\sqrt{L^2}-\sqrt{L^2}=0.$$
        Hence by the Hodge index theorem, we have $$\left(\gamma-\frac{L}{\sqrt{L^2}}\right)^2<0.$$ But $$\left(\gamma-\frac{L}{\sqrt{L^2}}\right)^2=\gamma^2+\frac{L^2}{L^2}-2\frac{L\cdot \gamma}{\sqrt{L^2}}=0.$$
        This is a contradiction. Hence $\gamma=\frac{[L]}{\sqrt{L^2}}$. Therefore, $\gamma_n\to \frac{L}{\sqrt{L^2}}$. The same argument applies to any convergent subsequence. Since $\gamma_n$ is bounded and every accumulation point equals $\frac{L}{\sqrt{L^2}}$, the whole sequence converges to $\frac{[L]}{\sqrt{L^2}}$.
    \end{proof}

    \begin{proposition}\label{prop:L^2 differentiable}
        Let $L$ be an ample line bundle on a surface $X$ with $\ex(L)=\sqrt{L^2}$.
Then  $\ex$ is differentiable at $L$ if and only if every curve $C$ computing $\ex(L)$ satisfies $$\frac{[C]}{{\rm mult}_x C}=\frac{[L]}{\sqrt{L^2}}.$$ 
In this case,
$$d\ex|_L(M)=\frac{L\cdot M}{\sqrt{L^2}}.$$
In particular, if no curve computes $\ex(L)$, then $\ex$ is differentiable at $L$.
    \end{proposition}

    \begin{proof}
      
    Suppose that every curve $C$ computing $\ex(L)$ satisfies $$\frac{[C]}{{\rm mult}_x C}=\frac{[L]}{\sqrt{L^2}}.$$

        For every nonzero $B$ of sufficiently small norm, both $L+B$ and $\frac{1}{2}L+B$ are ample. Choose an integral curve $C_{B}$ such that 
        \begin{equation}\label{eqn3}
            \ex(L+B)\le\frac{(L+B)\cdot C_{B}}{{\rm mult}_x C_{B}}\le \ex(L+B)+\norm{B}^2.
        \end{equation}
        First, we claim that $\frac{B\cdot C_B}{\mx C_B}\to 0$ as $\norm B  \to 0$. Indeed, $$\frac{(L+B)\cdot C_B}{{\rm mult}_x C_B}=\frac{1}{2}\frac{L\cdot C_B}{{\rm mult}_x C_B}+\frac{(\frac{1}{2}L+B)\cdot C_B}{{\rm mult}_x C_B}\ge \frac{1}{2}\frac{L\cdot C_B}{{\rm mult}_x C_B}.$$
        The last inequality follows from the fact that $\frac{1}{2}L+B$ is ample.
        By the choice of $C_B$ we can bound $\frac{(L+B)\cdot C_B}{{\rm mult}_x C_B}$. Hence, since $\frac{1}{2}L$ is ample, using \cref{lem:bdd numerical class} we can see that $\frac{ C_B}{{\rm mult}_x C_B}$ is bounded in $N_1(X)_{\mb R}$. Hence we have the claim as $|\frac{B\cdot C_B}{{\rm mult}_x C_B}|\le \lambda \norm{B}\norm{\frac{C_B}{{\rm mult}_x C_B}}$ for some constant $\lambda>0$ independent of $B$.
        Letting  $\norm{B}\to 0$ in \eqref{eqn3} by continuity of $\ex$ and $\|\cdot\|$, we get
        $$\frac{L\cdot C_{B}}{{\rm mult}_x C_{B}} \to \ex(L)=\sqrt{L^2}.$$
        From \cref{lem:seq argument} we get $$\frac{C_{B}}{{\rm mult}_x C_{B}}\to \frac{L}{\sqrt{L^2}}.$$

        \begin{equation*}
            \begin{split}
                 \ex(L+B)&\ge\frac{(L+B)\cdot C_{B}}{{\rm mult}_x C_{B}}-\norm{B}^2\\
                &=\frac{L\cdot C_{B}}{{\rm mult}_x C_{B}}+\frac{B\cdot C_{B}}{{\rm mult}_x C_{B}}-\norm{B}^2\\
                &\ge \sqrt{L^2}+\frac{L\cdot B}{\sqrt{L^2}}-\frac{L\cdot B}{\sqrt{L^2}}+\frac{B\cdot C_{B}}{{\rm mult}_x C_{B}}-\norm{B}^2\\
                &= \sqrt{L^2}+\frac{L\cdot B}{\sqrt{L^2}}+B\cdot\left(\frac{ C_{B}}{{\rm mult}_x C_{B}}-\frac{L}{\sqrt{L^2}}\right)-\norm{B}^2.
            \end{split}
        \end{equation*}
        
        We have
        
        $$\ex(L+B)\le\sqrt{(L+B)^2}\
                =\sqrt{L^2}+\frac{L\cdot B}{\sqrt{L^2}}+O(\norm{B}^2) \text{ (Taylor expansion). }$$

        This gives
        $$\lim_{\norm{B}\to 0}\frac{\ex(L+B)-\ex(L)-\frac{L\cdot B }{\sqrt{L^2}}}{\norm{B}}=0.$$

        Thus $\ex$ is differentiable at $L$.

        Conversely, suppose $\ex$ is differentiable at $L$. If no curve computes $\ex(L)$, the required condition is vacuous. Let $C$ be an integral curve such that $$\sqrt{L^2}=\ex(L)=\frac{L\cdot C}{{\rm mult}_x C}.$$

        From the definition of $\varepsilon(L,x)$ and $\ex(L)=\sqrt{L^2}$, the function
        $f(A)=\sqrt{A^2}-\ex(A)$ has a local minimum at $L$ (note that $\ex(A)\le \sqrt{A^2}$). Then we have $df_L=0$. Hence 
        \begin{equation}\label{eqn1}
            d\ex|_L(M)=\frac{L\cdot M}{\sqrt{L^2}}\text{ for all $M\in N_1(X)_{\mb R}$.}
        \end{equation}
        If $C$ computes $\ex(L)$, the function $\frac{A\cdot C}{{\rm mult}_x C}-\ex(A)$ has a local minimum at $L$. Hence we have

        \begin{equation}\label{eqn2}
            d\ex|_L(M)=\frac{M\cdot C}{{\rm mult}_x C}\text{ for all $M\in N_1(X)_{\mb R}$.}
        \end{equation}
        From \eqref{eqn1} and \eqref{eqn2} we can conclude
$$\frac{[C]}{{\rm mult}_x C}=\frac{[L]}{\sqrt{L^2}}.$$ This completes the proof.
    \end{proof}

    \begin{theorem}\label{thm: N_x nowhere dense}
        The set $N_x$ is nowhere dense in ${\rm Amp}_{\mb R}(X)$.
    \end{theorem}
    \begin{proof}
    To prove that $N_x$ is nowhere dense, it suffices to show that every non-empty open set contains a non-empty open subset on which $\ex$ is differentiable everywhere. Let $O$ be a non-empty open subset of ${\rm Amp}_{\mb R}(X)$. Suppose $O$ contains a line bundle $L$ such that $\ex(L)<\sqrt{L^2}$.
       By \cite[Proposition 1.3]{Sch} there exist curves $C_1,C_2,\dots,C_k$ and a neighborhood $U\subset {\rm Nef}_{\mb R}(X)$ of $L$ such that 
        \begin{equation}\label{eqn:local sc submaximal}
            \ex(L')=\min_{1\le i\le k} \frac{L'\cdot C_i}{{\rm mult}_x C_i}
        \end{equation}
        for all $L'\in U$. The $C_i$ are Seshadri curves of $L$ at $x$. We can choose $U$ to be inside $O$. From \cref{lem:seshadri curve class different}, we can see that $L\mapsto \frac{L\cdot C_i}{\mx C_i}$ defines a different linear map for each $i$. Consider an $L_0$ in $U$ such that $L_0$ is outside all  hyperplanes $H_{ij}$ for all $1\le i<j\le k$ where  $$H_{ij}=\left\{M\middle|\frac{M\cdot C_i}{\mx C_i}=\frac{M\cdot C_j}{\mx C_j}\right\}.$$
        Such an $L_0$ exists because the linear functionals are distinct by \cref{lem:seshadri curve class different}. Let $C$ be the curve in $\{C_1,C_2,\dots,C_k\}$ such that $\ex(L_0)=\frac{L_0\cdot C}{\mx C}$.
        Note that $C$ is unique by the choice of $L_0$. By the continuity of linear functionals ($\frac{\_\cdot C_i}{\mx C_i}$) and the choice of $L_0$ we can find an open set $V$ in $U$ such that $\frac{M\cdot C}{\mx C}<\frac{M \cdot C_i}{\mx C_i}$ for $C_i\ne C$. 
        Hence, by \eqref{eqn:local sc submaximal}, we have $\ex(M)=\frac{M\cdot C}{\mx C}$ for all $M\in V$. Since $\ex$ is  linear   on the open neighborhood $V$,  it is differentiable throughout $V$. Hence given any open set $U$ containing $L$, we can find an open set $V$ inside $U$ such that $\ex$ is differentiable throughout $V$.

        Now suppose $\ex(L)=\sqrt{L^2}$ for all $L\in O$. Since $L^2>0$ for ample bundles, $\ex$ is smooth on $O$. In both cases, $O$ contains a non-empty open subset $V$ disjoint from $N_x$. Since $V$ is open, it is also disjoint from $\overline{N_x}$. Therefore $N_x$ is nowhere dense in the ample cone.
    \end{proof}

        \begin{proposition}
            Let $X$ be a surface and $L$ be an ample line bundle on $X$ such that $\ex(L)=\sqrt{L^2}$. Suppose that $\ex$ is differentiable at $L$ and $K_X\cdot L+\sqrt{L^2}<0$. Then there exists an open neighborhood $U\subset {\rm Amp}_{\mb R}(X)$ of $L$ such that:
            $$\ex(L')=\sqrt{L'^2} \text{ for all } L'\in U.$$
        \end{proposition}
        \begin{proof}
            Suppose, to the contrary, that there exists a sequence of submaximal line bundles $L_n$ such that $L_n \to L$. Choose an integral curve $C_n$  such that
            $$\ex(L_n)=\frac{L_n\cdot C_n}{\mx C_n}.$$
            Since $\ex$ is continuous, we have $\ex(L_n)\to \ex(L)$. 

            We claim that $\frac{[C_n]}{\mx C_n}$ is bounded. Since $L_n-\frac{1}{2}L$ converges to the ample class $\frac12L$, the class $L_n-\frac{1}{2}L$ is ample for all sufficiently large $n$. Hence for large $n$ we have
            $$0\le\frac{1}{2}L\cdot\frac{C_n}{\mx C_n}
\le L_n\cdot\frac{C_n}{\mx C_n}
=\ex(L_n).$$
Since the sequence $(\varepsilon_x(L_n))$ converges and is
therefore bounded, we conclude that
$\frac{L}{2}\cdot\frac{C_n}{{\rm mult}_xC_n}$
is bounded. Now use \cref{lem:bdd numerical class} to see that $\frac{[C_n]}{\mx C_n}$ is bounded.

Now,
            
                $$\frac{L\cdot C_n}{\mx C_n}=\frac{L_n\cdot C_n}{\mx C_n}+\frac{(L-L_n)\cdot C_n}{\mx C_n}=\ex(L_n)+\frac{(L-L_n)\cdot C_n}{\mx C_n}\to \ex(L)=\sqrt{L^2}.$$
            The last convergence holds as $L_n \to L$ and $\frac{[C_n]}{\mx C_n}$ are bounded.
            From \cref{prop:L^2 differentiable} we can see that $C_n$ is a sequence of integral curves satisfying the hypotheses of \cref{lem:seq argument} and hence $\frac{[C_n]}{\mx C_n}$ converges to $\frac{[L]}{\sqrt{L^2}}$.
            Let $\widetilde{X}$ denote the blow-up of $X$ at $x$. It is known that if $\widetilde{C_n}$ denotes the strict transform of $C_n$ on $\widetilde{X}$, $\widetilde{C_n}^2<0$. 
Thus we have $$C_n^2-(\mx C_n)^2\le -1.$$
            This gives:
            $$\frac{1}{(\mx C_n)^2}\le 1-\frac{C_n^2}{(\mx C_n)^2}\to 0.$$
            The last convergence holds because $\frac{[C_n]}{\mx C_n}$ converges to $\frac{[L]}{\sqrt{L^2}}$ and hence $\frac{C_n^2}{(\mx C_n)^2}\to 1.$
            Hence we can conclude that $\mx C_n \to \infty$.

            By the adjunction formula, 
            $$K_X\cdot C_n+\mx C_n=K_{\widetilde{X}}\cdot \widetilde{C_n}\ge -2-\widetilde{C_n}^2\ge -1.$$
            Thus,
            $$\frac{K_X\cdot C_n}{\mx C_n}+1\ge \frac{-1}{\mx C_n}.$$
           Taking limits, we obtain
            $$\frac{K_X\cdot L}{\sqrt{L^2}}+1\ge 0.$$
            This contradicts our assumption. This completes the proof.
        \end{proof}

    \section{Maximality of Line Bundles on General Blow-Ups of $\mathbb{P}^2$}\label{sec:Submaximality}
	
	In this section, we study Seshadri functions on blow-ups of $\mb P^2$. We first construct an open submaximality region for arbitrary configurations of distinct points, and then obtain open maximality regions under suitable SHGH assumptions. From the continuity of the Seshadri function, it is clear that the set of all line bundles $L$ with $\varepsilon_x(L)<\sqrt{L^2}$ is open. We start by explicitly exhibiting an open set of such line bundles.
	\begin{lemma}\label{lessopen}
		For $r\ge 1$, let $X_r$ be a  blow-up of r distinct points in $\mb P ^2$ and let $L=eH-n_1E_1-\dots-n_rE_r$ be a line bundle in the ample cone of $X_r$. If $ 2n_1e>\left(\sum_{i=1}^{r} n_i^2\right) + n_1^2$, then $\varepsilon(X_r,L,x) < \sqrt{L^2}$ for all $x \in X_r$.
		\begin{proof}
			As $ 2n_1e>\left(\sum_{i=1}^{r} n_i^2\right) + n_1^2$ we have $(e-n_1)^2<\left( e^2-\sum_{i=1}^{r} n_i^2\right)$.  This implies $L\cdot (H-E_1)<\sqrt{L^2}$. Now, since dim $|H-E_1|>0$, there is a curve passing through any $x\in X_r$. So we get $\varepsilon(X_r,L,x) < \sqrt{L^2}$.
		\end{proof}
	\end{lemma}
	
	\begin{proposition}
		For $r\ge 1$, let $X_r$ be a  blow-up of r distinct points in $\mb P ^2$. There exists a non-empty open set $U \subset \rm{Nef}_{\mathbb{R}}(X_r)$ such that $\varepsilon(X_r,L,x) < \sqrt{L^2}$ for all $ L \in U$ and for all $x\in X_r$.
		\begin{proof}
			Let $$S=\left\{ L\in \text{\rm Nef}_{\mathbb{R}}(X_r)~~|~~ L=eH-n_1E_1-\dots-n_rE_r \text{ with } 2n_1e>\left(\sum_{i=1}^{r} n_i^2\right) + n_1^2\right\}.$$ By Lemma \ref{lessopen}, $U=S\cap \text{Amp}_{\mathbb{R}}(X_r)$ is an open set in $\text{\rm Nef}_{\mathbb{R}}(X_r)$ such that $\varepsilon(X_r,L,x) < \sqrt{L^2}$ for all $ L \in U$. It remains to show that $U$ is non-empty. Let $L=dH-n_1E_1-\dots-n_rE_r$ be an ample line bundle. It is clear that $eH-n_1E_1-\dots-n_rE_r$ is ample $\forall e>d$. Thus we can increase $e$ so that $2n_1e>\sum_{i=1}^{r} n_i^2+n_1^2$ is satisfied. Hence, $S \cap \text{Amp}_{\mathbb{R}}(X_r)$ is non-empty.
		\end{proof}
	\end{proposition}
	For a surface $X$ and $x\in X$, define the set $$W_x=\left\{ L\in \text{\rm Nef}_{\mathbb{R}}(X)~~| ~~\varepsilon(L,x)=\sqrt{L^2}\right\}.$$ When $X$ is an abelian surface with Picard rank $\rho(X)\ge 2$, $W_x$ is nowhere dense \cite[Lemma 4.1]{Sch}. We can ask a similar question for $X_r$. Let $W_x$ denote the same set for $X_r$. The following propositions conjecturally show  that $W_x$ is \textit{not} nowhere dense in ${\rm  Nef}_{\mb R}(X_r)$.
	
	For a line bundle $F$ in standard form, it is classical that $F\cdot C\geq 0$ for all $(-1)$-curves $C$. Although this is well known, we include a proof for completeness and for the convenience of the reader.
	
	\begin{lemma}\label{stdform}
    For $r\ge 3$,
		let $L=eH-n_1E_1-\dots-n_rE_r$ be a real line bundle in standard form. Then $L\cdot  C\geq 0$ for every $(-1)$-curve $C$ .
		\begin{proof}
			If $C$ is any of the exceptional divisors $E_1,\dots,E_r$, then $L\cdot C=L\cdot E_i=n_i\geq 0$. Now assume that  $C$  is a $(-1)$-curve and $C\neq E_i$ for $1\leq i\leq r$. We can write $L$ as a non-negative real combination of 
			$$H, H_1=H-E_1, H_2=2H-E_1-E_2, H_i=3H-E_1-E_2-\dots-E_i, ~~\forall 3\leq i\leq r $$
			as follows: 
			$$L=(e-(n_1+n_2+n_3))H+(n_1-n_2)H_1+(n_2-n_3)H_2+\sum_{i=3}^{r-1}(n_i-n_{i+1})H_i+n_rH_r.$$
			Note that for $3 \leq i \leq r$, $H_i+K_{X_r}=E_{i+1}+\dots + E_r$. As $C\neq E_i$, we have $C\cdot E_i\geq 0$. So $H_i\cdot C \geq -K_{X_r}\cdot C$. Note that $H,H_1,H_2$ are nef and since $C$ is a $(-1)$-curve, we have $-K\cdot  C=1$. So $L\cdot  C\geq 0$. 
		\end{proof}
	\end{lemma}
	
	\begin{corollary}\label{standard}
		For $r\ge 3$,
		let $L=eH-n_1E_1-\dots-n_rE_r$ be a real line bundle on $X_r$ with $n_i\geq 0$ and $e \geq n_i+n_j+n_k$ for all distinct $i,j,k \in \{1,2,\dots,r\}$. Then $L\cdot  C\geq 0$ for every $(-1)$-curve $C$ .
	\end{corollary}
	\begin{proof}
		We can convert such an $L$ into standard form using $\sigma_i\in W_r$. Also note that a $(-1)$-curve is again a $(-1)$-curve under the action of the Weyl group $W_r$. Then the result follows immediately from the previous lemma.
	\end{proof}
	
	From now on, $X_r$ will denote a  blow-up of $r$ very general points in $\mb P^2$. We will now prove some results about the Seshadri function for a very general point $x\in X_r$ assuming the Weak SHGH conjecture (Conjecture \ref{WSHGH}).
	
	\begin{theorem}\label{open}
		For $r\geq 9$, assume that \cref{WSHGH} holds for $X_{r+1}$.
Then, for a very general point $x\in X_r$, there exists a non-empty
open subset $U_x\subset {\rm Amp}_{\mb R}(X_r)$ such that
$$  \varepsilon(X_r,L,x)=\sqrt{L^2}
    \qquad\text{for every }L\in U_x.$$
In particular, the Seshadri function $\varepsilon_x$ is not
piecewise linear. Moreover, $U_x$ contains the numerical class
of an integral ample line bundle whose Seshadri constant at $x$
is irrational.
		\begin{proof}
			Since $L\mapsto \sqrt{L^2}$ is not a piece-wise linear function on $\text{\rm Nef}_{\mathbb{R}}(X_r)$, it is enough to construct an open set $U_x$ of $\text{\rm Nef}_{\mathbb{R}}(X_r)$ such that $\varepsilon(X_r,L,x)=\sqrt{L^2}$ for all line bundles $L\in U_x$. 
			
			Consider the following subset of $\text{Pic}_{\mathbb{R}}(X_r)$,
			$$V_x=\{L=eH-n_1E_1-\dots-n_rE_r~~|~~\sqrt{L^2}>n_i>0 \forall 1\le i \le r \text{ and } e>\sqrt{L^2}+n_j+n_k~~ \forall 1\leq j<k\leq r  \}.$$ Clearly $V_x$ is open in $\text{Pic}_{\mathbb{R}}(X_r)$. Now we will prove that $U_x=V_x\cap \text{Amp}_{\mathbb{R}}(X_r)$ has the desired property. Consider $L=eH-n_1E_1-\dots-n_rE_r \in U_x$ with $e,n_1,\dots,n_r \in \mathbb{Q}$. If $\varepsilon(L,x)<\sqrt{L^2}$, then $\varepsilon(nL,x)<\sqrt{(nL)^2}$ for any positive integer $n$. Choose $n\in \mathbb{N}$ such that $nL$ is an integral line bundle. So there is a curve $C$ containing $x$ such that $\frac{nL\cdot  C}{\text{mult}_xC}=\varepsilon(nL,x)$. The same $C$ will also serve as a Seshadri curve for $L$.
			
			Now we argue as in the proof of \cite[Proposition 2.1]{Han1}. Let $\pi_x:\Xr\rightarrow X_r$ be the  blow-up of $X_r$ at $x$. As $C$ is a Seshadri curve, $\widetilde{C}^2<0$, where $\widetilde C$ is the strict transform of $C$. Since we assumed that Conjecture \ref{WSHGH} is true, $\widetilde{C}$ is a $(-1)$-curve. Then \begin{equation}\label{less}
				0=(\pi_x^\ast(L)-\varepsilon(L,x) E_x)\cdot  \widetilde{C}>(\pi_x^\ast(L)-\left(\sqrt{L^2}\right) E_x)\cdot  \widetilde{C}.
			\end{equation} 
			By the choice of $L$, we have $\sqrt{L^2}>n_i$ and $e>\sqrt{L^2}+n_i+n_j$, $\forall 1\leq i\neq j \leq r$. Therefore by Corollary \ref{standard}, we have $(\pi_x^\ast(L)-\sqrt{L^2} E_x)\cdot \widetilde{C} \geq 0$ as $\widetilde{C}$ is a $(-1)$-curve. This is a contradiction to \eqref{less}. Hence $\varepsilon(L,x)=\sqrt{L^2}$ for every $L\in U_x\cap \text{Pic}_\mathbb{Q}(X_r)$. So the continuous function $\varepsilon_x(L)$ is equal to $\sqrt{L^2}$ for all $L$ in the dense set $U_x\cap \text{Pic}_\mathbb{Q}(X_r)$. Hence the Seshadri constant is maximal  for all line bundles in $U_x$. 
			
			It remains to prove that $U_x$ is non-empty. For  $r>12$ and  $r\neq 15,16$, it is easy to show that there exists $e$ such that $4e-4<r<e^2-1$. Consider $L=eH-E_1-\dots-E_r$ for such an $e$. It is easy to see that the Nagata conjecture implies that a uniform line bundle $L$ is ample if $L^2>0$. Note that Conjecture \ref{WSHGH} implies the Nagata conjecture. Since $\sqrt{L^2}=\sqrt{e^2-r}>0$, $L$ is an ample line bundle. To show that $L\in U_x$, it suffices to prove $e>\sqrt{e^2-r}+2$, which is equivalent to $r>4e-4$. So for $r>12$, ~  $r\neq 15,16$ the theorem is true. Furthermore, choose $e=\lceil r/4\rceil$ and let $L=eH-E_1-\dots-E_r$. For $r\geq 17$, we have $e\geq 5$ and $4e-3\leq r\leq 4e$,
so
$$(e-3)^2<e^2-4e\leq L^2\leq e^2-4e+3<(e-2)^2.$$
In particular, $L^2>1$ and $L^2$ is not a perfect square, and
$4e-4<r<e^2-1$. For $r=13,14$, the choice $e=4$
gives $L^2=3,2$, respectively.

            For $r \in \{9,10,11,12,15,16\}$, we  explicitly give an element in $U_x$ below:\\
			$r=9~$  : $~16H-5E_1-\dots-5E_r$,\\
			$r=10$ : $10H-3E_1-\dots-3E_r$,\\
			$r=11$ : $7H-2E_1-\dots-2E_r$,\\
			$r=12$ : $11H-3E_1-\dots-3E_r$,\\
			$r=15$ : $9H-2E_1-\dots-2E_r$,\\
			$r=16$ : $9H-2E_1-\dots-2E_r$.
			
			It is easy to check that the above line bundles satisfy the conditions to be in $U_x$. For example, consider the case $r=9$. Let $L=16H-5E_1-\dots-5E_r$. Then $L$ is in $U_x$ as $\sqrt{31}=\sqrt{L^2}>n_i=5$ and $16=e>\sqrt{L^2}+n_j+n_k=\sqrt{31}+5+5$. Other cases are similar.
		\end{proof}
	\end{theorem}
	
	Now consider $W_x'=\{ L \in \text{\rm Nef}_{\mathbb{R}}(X_r)~~|~~ \varepsilon(X_r,L,x) \notin \mathbb{Q}\}$. We ask: are there infinitely many line bundles in $W'_x$ having pairwise numerically non-proportional classes?
	\begin{corollary}\label{irrSC}
		For $r\geq9$, assume that Conjecture \ref{WSHGH} for $X_{r+1}$ is true. Then there exist infinitely many integral ample line bundles on $X_r$ whose numerical classes are pairwise non-proportional and whose Seshadri constants at a very general point $x\in X_r$ are irrational.
		\begin{proof}
			Let $L$ denote the ample line bundle that we constructed in Theorem \ref{open}, with an irrational Seshadri constant. Thus, for sufficiently large $n$, $L+\frac{1}{n}H$ is ample and belongs to $U_x$. Hence the ample line bundle $L_n=nL+H$ has Seshadri constant $\sqrt{L_n^2}$. 
            Suppose, to the contrary, that $\sqrt{L_n^2}$ is irrational for at most finitely many $n$. Hence, for all sufficiently large $n$, the number
$$
\sqrt{L_n^2}=\sqrt{n^2L^2+2nL\cdot H+1}
$$
is a positive integer.
            The sequence
            $$\sqrt{L_{n+1}^2}-\sqrt{L_n^2}=\frac{(2n+1)L^2+2L\cdot H}
{\sqrt{L_{n+1}^2}+\sqrt{L_n^2}}=\frac{(2+1/n)L^2+2(L\cdot H)/n}
{\sqrt{L_{n+1}^2/n^2}+\sqrt{L_n^2/n^2}}$$
converges to $\sqrt{L^2}$ as $n\to \infty$. This is a contradiction, as the sequence of integers converges to an irrational number.

Suppose $L_n=\lambda L_m$, so that $nL+H=\lambda(mL+H)$.
Write $L=eH-m_1E_1-\dots-m_rE_r$, where $m_i>0$ by ampleness.
Comparing the coefficients of $E_i$ gives $n=\lambda m$,
while comparing the coefficients of $H$ gives
$ne+1=\lambda(me+1)=ne+\lambda$.
Hence, $\lambda=1$ and $n=m$. Therefore, the classes $L_n$
are pairwise non-proportional.
		\end{proof}
	\end{corollary}
	The existence of infinitely many ample line bundles with irrational Seshadri constants under the Weak SHGH conjecture was already observed in \cite[Remark~2.5]{Han1}; \cref{irrSC} gives an explicit construction of infinitely many pairwise numerically non-proportional such classes.
	
	Now we exhibit some results regarding the Seshadri function for an arbitrary $x\in X_r$ assuming the Strong SHGH conjecture (Conjecture \ref{SSHGH}).
	
	\begin{lemma}\label{std1}
		Assume that Conjecture \ref{SSHGH} is true and $r\ge 3$. For $x\in X_r$, let $\Xr$ be the  blow-up of $X_r$ at $x$ and let $E_x$ denote the corresponding exceptional divisor. Let $L=eH-n_1E_1-\dots-n_rE_r-n_xE_x\in \rm{Pic}_\mathbb{R}(\Xr)$ be a line bundle on $\Xr$ with $e \geq n_i+n_j+4n_k$, $n_i \geq n_x$ and $n_x\geq 0$ for all distinct $i,j,k \in \{1,\dots ,r\}$. Then $L\cdot C\geq 0$ for all irreducible, reduced curves $C$ with $C^2 < 0$.
		\begin{proof}
			If $C=E_i$ for some $i\in \{1,2,\dots ,r\}$, then $L\cdot C=L\cdot E_i=n_i\geq 0$. If $C=E_i-E_x$ for some $i\in \{1,2,\dots ,r\}$, then $L\cdot C=L\cdot (E_i-E_x)=n_i-n_x\geq 0$. Also $L\cdot E_x=n_x\ge 0$. Now assume that $C$ is any curve other than $E_i$ or $E_i-E_x$ or $E_x$, with $C^2<0$. 
		
			Choose a permutation $\sigma\in S_r$ such that
$n_{\sigma(1)}\geq\cdots\geq n_{\sigma(r)}$.
Set $n'_i=n_{\sigma(i)}$ and $E'_i=E_{\sigma(i)}$
for $1\leq i\leq r$. Now we can write $L'=eH-n'_1E'_1-\cdots-n'_rE'_r-n_xE_x$ as a non-negative real combination of $H, H_1=H-E_1', H_2=2H-E_1'-E_2', H_i=3H-E_1'-\dots-E_i'$, $\forall 3\leq i\leq r$, and $H_x=-K_{\Xr}$ as below:
			$$L'=(e-(n_1'+n_2'+n_3'))H+(n_1'-n_2')H_1+(n_2'-n_3')H_2+ \sum_{i=3}^{r-1}(n'_i-n'_{i+1})H_i+(n_r'-n_x)H_r+n_xH_x. $$
			Note that $L'=L$. As we explained in the proof of Lemma \ref{stdform}, $H_i\cdot  C\geq -K\cdot  C$ for $3\leq i \leq r$ and $H,H_1,H_2$ are nef. Also note that $H\cdot C \geq 1$. By \cite[Lemma 5.1]{Cyr}, $-K_{\Xr}\cdot C\geq -3$. So we have $H_i\cdot C \geq -3$ for $3\leq i \leq r$. Therefore,
			\begin{equation*}
				\begin{split}
					L\cdot C&=L'\cdot C\\
&=(e-n'_1-n'_2-n'_3)H\cdot C
 +(n'_1-n'_2)H_1\cdot C
 +(n'_2-n'_3)H_2\cdot C\\
&\quad+\sum_{i=3}^{r-1}(n'_i-n'_{i+1})H_i\cdot C
 +(n'_r-n_x)H_r\cdot C+n_xH_x\cdot C\\
&\ge e-n'_1-n'_2-n'_3
 -3\sum_{i=3}^{r-1}(n'_i-n'_{i+1})
 -3(n'_r-n_x)-3n_x\\
&=e-n'_1-n'_2-4n'_3\\
&\ge0.
				\end{split}				
			\end{equation*}
			The last inequality follows from the hypotheses.
		\end{proof}
	\end{lemma}
	
	\begin{theorem}
		For $r\geq 36$, assume that \cref{SSHGH} holds for $X_r$.
Then there exists a non-empty open subset
$U\subset {\rm Amp}_{\mb R}(X_r)$ such that
$$\varepsilon(X_r,L,x)=\sqrt{L^2}
    \qquad
    \text{for every }L\in U\text{ and every }x\in X_r.$$
		\begin{proof}
			Consider the following subset of $\text{Pic}_{\mathbb{R}}(X_r)$,
			$$V=\{L=eH-n_1E_1-\dots-n_rE_r~~ | ~~L^2>0, \sqrt{L^2}<n_j, e>n_j+n_k+4n_l \text{ for distinct } j,k,l \in \{1,\dots ,r \}  \}.$$ Clearly $V$ is open in $\text{Pic}_{\mathbb{R}}(X_r)$. We will prove that $U=V\cap \text{Amp}_{\mathbb{R}}(X_r)$ has the desired property.
			
			First, we  show that $U$ is non-empty. For $r\geq 36$, consider $e\in \mathbb{Q}_+$ such that $$r(r-1)^2<e^2<(r+1)(r-1)^2\text{ and }L=eH-(r-1)E_1-\dots-(r-1)E_r.$$
			
			As before, we can show that $L$ is ample, since $L$ is uniform and $L^2>0$ and we have $e>\sqrt{r}(r-1)\geq 6(r-1)$. To show that $L\in U$, it remains to prove that $\sqrt{e^2-r(r-1)^2}<r-1$, which is equivalent to $e^2<(r+1)(r-1)^2$. This is true by the choice of $e$. So $U$ is non-empty.
			
			Fix an arbitrary point $x\in X_r$. Consider $L=eH-n_1E_1-\dots-n_rE_r \in U$ with $e,n_1,\dots,n_r \in \mathbb{Q}_{+}$. Suppose that $\varepsilon(L,x)<\sqrt{L^2}$. The same is true for any positive integral multiple of $L$. Let $n$ be a positive integer such that $nL$ has integer coefficients. So there is a curve $C$ such that $\frac{nL\cdot  C}{\text{mult}_xC}=\varepsilon(nL,x)$. By homogeneity, $C$ also computes $\varepsilon(L,x)$.
			
			Let $\pi_x:\Xr\rightarrow X_r$ be the  blow-up of $X_r$ at $x$ and let $E_x$ denote the exceptional divisor. As $C$ is a Seshadri curve, we have $\widetilde{C}^2<0$, where $\widetilde{C}$ is the strict transform of $C$.
			Then
			\begin{equation}\label{less1}
				0=(\pi_x^\ast(L)-\varepsilon(L,x) E_x)\cdot  \widetilde{C}>(\pi_x^\ast(L)-\sqrt{L^2} E_x)\cdot  \widetilde{C}.
			\end{equation} 
			By the choice of $L$, $\sqrt{L^2}<n_i$ and $e>n_i+n_j+4n_k$ for distinct $i,j,k \in \{1,\dots ,r\}$. Therefore by Lemma \ref{std1}, we have $(\pi_x^\ast(L)-\sqrt{L^2} E_x)\cdot \widetilde{C} \geq 0$. This is a contradiction to \eqref{less1}. Hence for every $L\in U\cap \text{Pic}_\mathbb{Q}(X_r)$, $\varepsilon(L,x)=\sqrt{L^2}$. Since $U\cap \text{Pic}_\mathbb{Q}(X_r)$ is dense in $U$, the Seshadri constant is maximal for all line bundles in $U$. 
		\end{proof}
	\end{theorem}

    \bibliographystyle{halpha}
	\bibliography{SC.bib}
\end{document}